\documentclass{amsproc}

\usepackage{}

\newtheorem{theorem}{Theorem}[section]
\newtheorem{lemma}[theorem]{Lemma}

\theoremstyle{definition}
\newtheorem{definition}[theorem]{Definition}
\newtheorem{example}[theorem]{Example}

\theoremstyle{remark}

\numberwithin{equation}{section}
\DeclareMathOperator{\Res}{Res}
\def\@adminfootnotes{
  \let\@makefnmark\relax  \let\@thefnmark\relax
  \ifx\@empty\@date\else \@footnotetext{\@setdate}\fi
  \ifx\@empty\@subjclass\else \@footnotetext{\@setsubjclass}\fi
  \ifx\@empty\@keywords\else \@footnotetext{\@setkeywords}\fi
  \ifx\@empty\thankses\else \@footnotetext{
    \def\par{\let\par\@par}\@setthanks}
  \fi
}
\makeatother
\begin{document}

\title[$Q$-systems and Orthogonal Polynomials]{Generalizations of $Q$-systems and Orthogonal Polynomials from Representation Theory}


\author{Darlayne Addabbo and Maarten Bergvelt}
\address{Department of Mathematics, University of Illinois at Urbana-Champaign}
\curraddr{}
\email{addabbo2@illinois.edu, bergv@illinois.edu}
\thanks{Thanks to Rinat Kedem and Philippe Di Francesco for their helpful comments.}

\begin{abstract} We briefly describe what tau-functions in integrable
  systems are. We then define a collection of tau-functions given as
  matrix elements for the action of $\widehat{GL_2}$ on two-component
  Fermionic Fock space. These tau-functions are solutions to a
  discrete integrable system called a $Q$-system.

  We can prove that our tau-functions satisfy $Q$-system relations by
  applying the famous ``Desnanot-Jacobi identity'' or by using
  ``connection matrices", the latter of which gives rise to orthogonal
  polynomials. In this paper, we will provide the background
  information required for computing these tau-functions and obtaining
  the connection matrices and will then use the connection matrices to
  derive our difference relations and to find orthogonal polynomials.

  We generalize the above by considering tau-functions that are matrix
  elements for the action of $\widehat{GL_3}$ on three-component
  Fermionic Fock space, and discuss the new system of discrete
  equations that they satisfy. We will show how to use the connection
  matrices in this case to obtain ``multiple orthogonal polynomials of
  type II".
\end{abstract}
\date{\today}
\maketitle
\section{Introduction}
Integrable differential equations, such as the KdV equation,
\begin{equation}
u_t+u_{xxx}+6uu_{x}=0,
\end{equation} can be solved exactly by employing a change of
variables to rewrite the equations more simply in bilinear form. In
the case of the KdV equation, this change of variables is given by
\cite{MR2085332} 
\begin{equation}u=2(\ln \tau)_{xx}.
\end{equation} 
This method of changing variables is called Hirota's method and the
solutions of these differential equations under the change of
variables are referred to as ``tau-functions" (For more details on
using Hirota's method to find solutions to the KdV equation, as well
as many other examples, see \cite{MR2085332}.)

Interestingly, tau-functions are often equal to matrix elements for
representations of infinite dimensional Lie groups (see, for example,
\cite{MR1736222} and \cite{MR1021978}).

In this paper, we will discuss tau-functions that satisfy discrete
integrable equations. We will first define tau-functions that are
given as matrix elements for the action of $\widehat{GL_2}$ on
two-component Fermionic Fock space and will discuss how to show that
these tau-functions satisfy a $Q$-system.

More specifically, we will see that our $\widehat{GL_2}$ tau-functions
satisfy 
\begin{equation}
\label{eq:2}
\tau_{k}^{(\alpha)}\tau_{k-2}^{(\alpha+2)}=\tau_{k-1}^{(\alpha+2)}\tau_{k-1}^{(\alpha)}-(\tau_{k-1}^{(\alpha+1)})^2,
\end{equation} for $k\ge 0$ and $\alpha\in \mathbb{Z}$. By applying a
suitable change of variables, this can be shown to be equivalent to
the defining relations for the $A_{\infty/2}$ $Q$-system which is
discussed, for example, in \cite{MR2566162}. These difference
relations are found using ``connection matrices" (defined below) and
these connection matrices can also be used to obtain orthogonal
polynomials. 

$Q$-systems are discrete integrable systems that appear in various
places in mathematics, for example, as the relations satisfied by
characters of Kirillov-Reshetikhin modules (see \cite{MR906858},
\cite{MR1255302}) or as mutations in a cluster algebra (see
\cite{MR2452184}, \cite{MR2551179}).

Since $Q$-systems and orthogonal polynomials are already interesting,
it is natural to ask what sort of discrete relations are satisfied by
analogous tau-functions, given as matrix elements for the action of
$\widehat{GL_3}$ on three-component Fermionic Fock space and what sort
of orthogonal polynomials come from the corresponding connection
matrices. In the following, we will describe how to define these new
$\widehat{GL_3}$ tau-functions and how to use connection matrices to
show that they satisfy the following system of equations, for all
$k,\ell\ge 0$ and $\alpha,\beta\in \mathbb{Z},$
\begin{align*}
  \tag{1} (\tau_{k,\ell}^{(\alpha+1,\beta)})^2=\tau_{k,\ell}^{(\alpha,\beta)}\tau_{k,\ell}^{(\alpha+2,\beta)}+\tau_{k+1,\ell+1}^{(\alpha,\beta)}\tau_{k-1,\ell-1}^{(\alpha+2,\beta)}-\tau_{k+1,\ell}^{(\alpha,\beta)}\tau_{k-1,\ell}^{(\alpha+2,\beta)}\\
  \tag{2}\tau_{k,\ell}^{(\alpha+1,\beta)}\tau_{k,\ell-1}^{(\alpha+2,\beta)}=\tau_{k-1,\ell-1}^{(\alpha+2,\beta)}\tau_{k+1,\ell}^{(\alpha+1,\beta)}+\tau_{k,\ell-1}^{(\alpha+1,\beta)} \tau_{k,\ell}^{(\alpha+2,\beta)}\\
  \tag{3} (\tau_{k,\ell}^{(\alpha,\beta+1)})^2=\tau_{k,\ell}^{(\alpha,\beta)}\tau_{k,\ell}^{(\alpha,\beta+2)}-\tau_{k,\ell-1}^{(\alpha,\beta+2)}\tau_{k,\ell+1}^{(\alpha,\beta)}-\tau_{k+1,\ell}^{(\alpha,\beta+2)}\tau_{k-1,\ell}^{(\alpha,\beta)}\\
  \tag{4} \tau_{k-1,\ell}^{(\alpha,\beta+1)}\tau_{k,\ell+1}^{(\alpha,\beta)}=\tau_{k-1,\ell}^{(\alpha,\beta)}\tau_{k,\ell+1}^{(\alpha,\beta+1)}+\tau_{k-1,\ell+1}^{(\alpha,\beta)}\tau_{k,\ell}^{(\alpha,\beta+1)}.
\end{align*}
We hope, similarly to $Q$-systems, that our new system of equations
will also have connections to other areas of mathematics. We will
briefly discuss progress we have made in analyzing this new system of
equations.

Applying restrictions to the connection matrices in the
$\widehat{GL_3}$ case, we find an analogous collection of orthogonal
polynomials, which we will discuss. In our future work, we hope to
investigate more general situations, obtained by dropping these
restrictions.

(See \cite{paperinprogress} for more details on the computations of
our tau-functions and the difference relations that they
satisfy. Orthogonal polynomials, however are not discussed there.)

\section{Calculating $\widehat{GL_2}$ Tau-Functions on Two-Component Fermionic Fock Space}
Before we define our $\widehat{GL_2}$ tau-functions, we first define
two-component Fermionic Fock space, $F^{(2)}$ and describe the action
of $\widehat{GL_2}$ on this space.  Here, we will omit most technical
details. For more information, we refer the reader to \cite{MR1099256}
and \cite{paperinprogress}. In particular, all omitted details of the
following background on Fermionic Fock space and the associated action
of $\widehat{gl_2}$ can be found in \cite{paperinprogress}.

Consider the vector space
$H^{(2)}:=\mathbb{C}^{2}\otimes \mathbb{C}[z,z^{-1}]$. A basis of this
space is given by elements, $e_{a}z^k$, $a=0,1$, $k\in \mathbb{Z}$,
where $e_0=
\begin{pmatrix}
  1\\
  0
\end{pmatrix}$ and $e_1=
\begin{pmatrix}
  0\\
  1
\end{pmatrix}$.  
$F^{(2)}$ is then spanned by vectors,
\[
w=w_0\wedge w_1\wedge w_2 \wedge\cdots,
\]
where $w_i\in H^{(2)}$ and the $w$ satisfy some restrictions that we
will now discuss.

Let the vacuum vector be 
\[
v_0:=
\begin{pmatrix}
  1\\
  0
\end{pmatrix}
\wedge
\begin{pmatrix}
  0\\
  1
\end{pmatrix}
\wedge
\begin{pmatrix}
  z\\
  0
\end{pmatrix}
\wedge
\begin{pmatrix}
  0\\
  z
\end{pmatrix}
\wedge\cdots
\in F^{(2)},
\]
and define operators, $e(e_a z^k)$ and $i(e_a z^k)$ (called exterior
and interior product operators, respectively), given by
$e(e_a z^k)w=e_a z^{k}\wedge w$ and $i(e_a z^k)w=\beta$ if
$w=e_a z^k \wedge \beta$. $F^{(2)}$ is the span of the vectors
obtained by acting on $v_0$ by finitely many exterior and interior
product operators. We can specify an order in which to act by these
exterior and interior product operators and define ``elementary
wedges" as those wedges obtained by acting on $v_0$ by monomials of
exterior and interior product operators, subject to this order. For
more information, see \cite{MR1099256} and \cite{paperinprogress}. The
elementary wedges are defined in such a way that there exists a unique
bilinear form on $F^{(2)}$, denoted $\langle v,w\rangle$ for
$v,w\in F^{(2)}$, for which the elementary wedges are an orthonormal
basis.

It is useful to introduce generating series, called fermion fields,
for the exterior and interior product operators.  Define
\[
\psi_{a}^{\pm}(w)=\displaystyle\sum_{k\in
  \mathbb{Z}}{_a}\psi^{\pm}_{(k)}w^{-k-1}, \quad a=0,1,
\] 
where
\[
{_a}\psi^{+}_{(k)}=e(e_a z^k)\text{ and }{_a}\psi^{-}_{(k)}=i(e_a z^{-k-1}).\]
We can use fermion fields to express the action of
$\widehat{gl_{2}}=gl_{2}\otimes \mathbb{C}[z,z^{-1}]\oplus\mathbb{C}c$
on $F^{(2)}$. Let $E_{ab}\in gl_2$ $(a,b=0,1)$ be the matrices such
that $E_{ab}e_c=\delta_{bc}e_a$ and let the current,
$E_{ab}(w)=\displaystyle\sum_{k\in \mathbb{Z}}E_{ab}z^k w^{-k-1},$ be
the generating series of elements in $\widehat{gl_2}$. When $a\ne b$,
the series acts on $F^{(2)}$ by
\begin{equation}
\label{eq:1}
E_{ab}(w)=\psi_{a}^{+}(w)\psi_{b}^{-}(w).
\end{equation} The action of $E_{ab}(w)$ in general requires using the
normal ordered product, but we do not discuss this here since it is
not needed in this paper 
(More details can be found in \cite{paperinprogress}.).

In addition to the action of the Lie algebra, $\widehat{gl_2}$,
$F^{(2)}$ also carries an action of the group $\widehat{GL_2}$, a
central extension of the loop group, $\widetilde{GL_2}$. In
particular, on $F^{(2)}$ we have the action of ``fermionic translation
operators" $Q_0,Q_1$ such that
\[
\pi(Q_0)=
\begin{bmatrix}
  z^{-1} & 0\\
  0 & 1
\end{bmatrix}
, \quad \pi(Q_1)=
\begin{bmatrix}
1 & 0\\
0 & z^{-1}
\end{bmatrix}
\] 
where $\pi$ denotes the projection from $\widehat{GL_2}$ to
$\widetilde{GL_2}$. We also define $T=Q_1 Q_0^{-1}$ such that
\[\pi(T)=
\begin{bmatrix}
  z & 0\\
  0 & z^{-1}
\end{bmatrix}.
\]
Let $g_C\in \widehat{GL_2}$ be such that
\[\pi(g_C)=
\begin{bmatrix}
  1 & 0\\
  C(z) & 1
\end{bmatrix},
\] 
where
$C(z)=\displaystyle\sum_{i \in
  \mathbb{Z}}\displaystyle\frac{c_i}{z^{i+1}}$,
$c_i \in \mathbb{C}$. (Actually, for our purposes, it will sometimes
be useful to take the $c_i$s to be formal variables.)

We define our (``unshifted") tau-functions to be
\[\tau_{k}=\langle T^kv_{0},g_C\cdot v_{0}\rangle.\]
We also need to define ``shifted'' tau-functions, corresponding to the
action of $g_C^{(\alpha)}$ on the vacuum vector, where
\[g_C^{(\alpha)}=Q_0^{\alpha}g_C Q_0^{-\alpha},\quad \alpha \in
\mathbb{Z}.\]
Denote by $C^{(\alpha)}$, the series
$C^{(\alpha)}(z)=\displaystyle\sum_{i \in
  \mathbb{Z}}\displaystyle\frac{c_{i+\alpha}}{z^{i+1}}$,
so that $\pi(g_C^{(\alpha)})=
\begin{bmatrix}
  1 & 0\\
  C^{(\alpha)}(z) & 1
\end{bmatrix}$. We can calculate these tau-functions by noting that
\[g_C^{(\alpha)}=\exp
\begin{bmatrix}
  0 & 0\\
  C^{(\alpha)}(z) & 0
\end{bmatrix}
=\exp(\Res_{w}(C^{(\alpha)}(w)E_{10}(w))),
\]
where the action of the current, $E_{10}(w)$ on $F^{(2)},$ is given by
\eqref{eq:1}.  We then have the following formulas for our
tau-functions, which are stated and proven in \cite{paperinprogress}:
\begin{theorem}
\begin{enumerate}
\item $\tau_{k}^{(\alpha)}=0$ for $k<0.$
\item $\tau_{0}^{(\alpha)}=1.$
\item When $k>0,$
\[\tau_{k}^{(\alpha)}=\frac{1}{k!}\Res_{\mathbf{w}}(\displaystyle \prod_{1\le i<j\le k}(w_i-w_j)^2\prod_{i=1}^{k}C^{(\alpha)}(w_i)),\]
where $\Res_\mathbf{w}$ denotes
$\Res_{w_1}(\Res_{w_2}\cdots(\Res_{w_k}\cdots)\cdots)$.
\end{enumerate}
Alternatively, when $k>0$ we can write
\[\tau_{k}^{(\alpha)}=\det\left[\begin{array}{cccc}
c_{\alpha} & c_{\alpha+1} & \cdots & c_{\alpha+k-1}\\
c_{\alpha+1} & c_{\alpha+2} & \cdots & c_{\alpha+k}\\
\vdots & \vdots & \cdots & \vdots\\
c_{\alpha+k-1} & c_{\alpha+k} & \cdots & c_{\alpha+2k-2}
\end{array}\right].\]
\end{theorem}

We note in particular that these tau-functions are determinants of
Hankel matrices and are thus especially well suited to applying the
famous Desnanot-Jacobi identity \cite{MR1718370}. We obtain
\begin{equation}
\label{eq:3}
\tau_{k}^{(\alpha)}\tau_{k-2}^{(\alpha+2)}=\tau_{k-1}^{(\alpha+2)}\tau_{k-1}^{(\alpha)}-(\tau_{k-1}^{(\alpha+1)})^2,
\end{equation} for $k\ge 0$ and $\alpha\in \mathbb{Z}.$

\section{Generalizing to $\widehat{GL_3}$}
We now generalize the above to the $\widehat{GL_3}$ case. Here, in the
same way that we did in the $\widehat{GL_2}$ case, we first define the
action of $\widehat{gl_3}$ on three-component Fermionic Fock space,
$F^{(3)}$.  We define
$H^{(3)}:=\mathbb{C}^3\otimes \mathbb{C}[z,z^{-1}]$, which has a basis
given by $e_a z^{k}$, where $a=0,1,2$, $k\in \mathbb{Z}$, and
$e_0=\begin{pmatrix}
  1\\
  0\\
  0\\
\end{pmatrix}$, $e_1=\begin{pmatrix}
0\\
1\\
0\\
\end{pmatrix}$, $e_2=\begin{pmatrix}
0\\
0\\
1\\
\end{pmatrix}.$

The vacuum vector is
\[v_0=\begin{pmatrix}
1\\
0\\
0\\
\end{pmatrix}\wedge \begin{pmatrix}
0\\
1\\
0\\
\end{pmatrix}\wedge \begin{pmatrix}
0\\
0\\
1\\
\end{pmatrix}\wedge \begin{pmatrix}
z\\
0\\
0\\
\end{pmatrix}\wedge \begin{pmatrix}
0\\
z\\
0\\
\end{pmatrix}\wedge \begin{pmatrix}
0\\
0\\
z\\
\end{pmatrix}\wedge \cdots\in F^{(3)}.\]
Similarly to the $\widehat{GL_2}$ case, in addition to the action of
$\widehat{gl_3}$ on $F^{(3)}$, we have an associated action of
$\widehat{GL_3}$, a central extension of the loop group,
$\widetilde{GL_3}$. We have as before, a projection map, $\pi$ from
$\widehat{GL_3}$ to $\widetilde{Gl_3}$. We define fermionic
translation operators, $Q_0,Q_1,$ and $Q_2$ such that
\[\pi(Q_0)=\begin{bmatrix}
z^{-1} & 0 & 0\\
0 & 1 & 0\\
0 & 0 & 1
\end{bmatrix},\quad \pi(Q_1)=\begin{bmatrix}
1 & 0 & 0\\
0 & z^{-1} & 0\\
0 & 0 & 1
\end{bmatrix},\quad \pi(Q_2)=\begin{bmatrix}
1 & 0 & 0\\
0 & 1 & 0\\
0 & 0 & z^{-1}
\end{bmatrix}.\] We also need operators, $T_1=Q_1 Q_0^{-1}$ and $T_2=Q_2 Q_1^{-1}$ such that 
\[\pi(T_1)=\left[\begin{array}{ccc}
z & 0 & 0\\
0 & z^{-1} & 0\\
0 & 0 & 1
\end{array}\right] \text{ and } \quad \pi(T_2)=\left[\begin{array}{ccc}
1 & 0 & 0\\
0 & z & 0\\
0 & 0 & z^{-1}
\end{array}\right].\] 
The fermion fields and action of $\widehat{gl_3}$ are defined in a way
completely analogous to how they were defined in our discussion of the
$\widehat{GL_2}$ case. As before, we omit any discussion of the normal
ordered product, since the details are not neccesary for this paper
and can be found in \cite{paperinprogress}.

We take a loop group element, $g_{C,D,E}\in \widehat{Gl_3}$, such
that \[\pi(g_{C,D,E})=\left[\begin{array}{ccc}
                              1 & 0 & 0\\
                              C(z) & 1 & 0\\
                              D(z) & E(z) & 1
                            \end{array}\right],
\] 
where $C(z)=\displaystyle\sum_{i\in\mathbb{Z}} \frac{c_i}{z^{i+1}}$,
$D(z)=\displaystyle\sum_{i\in \mathbb{Z}} \frac{d_i}{z^{i+1}}$, and
$E(z)=\displaystyle\sum_{i\in \mathbb{Z}} \frac{c_i}{z^{i+1}}$ and the
$c_i$s, $d_i$s, and $e_i$s are complex numbers or formal variables.

We define
\[\tau_{k,\ell}=\langle g_{C,D,E}\cdot
v_{0},T_{1}^kT_{2}^{\ell}v_{0}\rangle ,\]
and calculate the action of $g_{C,D,E}$ in the same way that we
calculated the action of our previous group element on the
two-component Fermionic Fock space, by expressing $g_{C,D,E}$ in terms
of fermion fields.

In order to obtain our difference relations we must, as in the
$\widehat{GL_2}$ case, introduce shifted tau-functions. Here, we have
two independent shifts. We define the
$\tau_{k,\ell}^{(\alpha,\beta)}$s to be the tau-functions
corresponding to the action of the group
element,
\[g_{C,D,E}^{(\alpha,\beta)}=Q_0^{-\alpha}Q_{1}^{-\beta}g_{C,D,E}Q_{1}^{\beta}Q_0^{\alpha}.\]
We comment that we do not need $Q_2$ to obtain all possible shifts
since
\[Q_{2} g_{C,D,E}Q_{2}^{-1}=Q_{0}^{-1}Q_{1}^{-1}g_{C,D,E}Q_{1}Q_{0},\]
so there really are only two independent shifts, as claimed.

The formula for our $\tau_{k,\ell}^{(\alpha,\beta)}$ functions is
then: (The following is stated and proven in \cite{paperinprogress}).
\begin{theorem}
\[\tau_{k,\ell}^{(\alpha,\beta)}=\displaystyle\sum_{n_c+n_d=k,n_e+n_d=\ell} c_{n_c,n_d,n_e}\] where 
\[c_{n_c,n_d,n_e}=\frac{1}{n_c!n_d!n_e!}\Res_{\mathbf{x}}\Res_{\mathbf{y}}\Res_{\mathbf{z}}\left(\displaystyle\prod_{i=1}^{n_c}C^{(\alpha-\beta)}(x_{i})\displaystyle\prod_{i=1}^{n_d}D^{(\alpha)}(y_{i})\displaystyle\prod_{i=1}^{n_e}E^{(\beta)}(z_{i})p_{n_c,n_d,n_e}\right)\]
where we use the same notation for residues as we did in the
$\widehat{GL_2}$ case and

\[p_{n_c,n_d,n_e}=\displaystyle
(-1)^{\displaystyle\frac{n_d(n_d+1)}{2}}\times\]
\[\frac{\displaystyle \prod_{1\le i<j\le n_c}(x_i-x_j)^2 \displaystyle \prod_{1\le i<j\le n_d}(y_i-y_j)^2 \displaystyle \prod_{1\le i<j \le n_e}(z_i-z_j)^2 \displaystyle \prod_{i=1}^{n_c}\prod_{j=1}^{n_d}(x_i-y_j)\displaystyle \prod_{i=1}^{n_d}\prod_{j=1}^{n_e}(y_i-z_j)}{\displaystyle \prod_{i=1}^{n_c}\prod_{j=1}^{n_e}(x_i-\underline{z_j})}.\]
\end{theorem} 
Here, the underline means that we expand
$\frac{1}{x_i-\underline{z_j}}$ in positive powers of $z_j$. We
comment that the formula given here is analogous to the formula for
our $\widehat{GL_2}$ tau-functions, but is much more complicated. In
particular, the denomenator appearing in our formula here means that
our tau-functions are, in general, infinite series of monomials in the
$c_i$s, $d_i$s, and $e_i$s.
\section{Connection Matrices and Zero Curvature Equations}
Although the Desnanot-Jacobi identity easily yields difference
relations for our $\widehat{GL_2}$ tau-functions, it is unclear how
one might use it to find difference relations for our $\widehat{GL_3}$
tau-functions. To find difference relations for our $\widehat{GL_3}$
tau-functions, we use ``connection matrices" and show that certain
zero curvature equations are satisfied. For this, we need the Birkhoff
factorization of $T_{2}^{-\ell}T_{1}^{-k}g_{C,D,E}^{(\alpha,\beta)}$,
which exists when $\tau^{(\alpha,\beta)}_{k,\ell}\ne 0$
(\cite{MR900587}). Here, we consider our $c_{i}$s, $d_{i}$s, and
$e_{i}$s to be formal variables and assume that the Birkhoff
factorization exists for all
$T_{2}^{-\ell}T_{1}^{-k}g_{C,D,E}^{(\alpha,\beta)}$ where $k$ and
$\ell$ are both nonnegative and $\alpha,\beta\in \mathbb{Z}$.

\subsection{$\widehat{GL_2}$ Case} In this section, we will use zero
curvature equations (\cite{MR905674}) to rederive our difference
relations for the $\widehat{GL_2}$ case. The $\widehat{GL_3}$ case is
similar, the details of which are included in \cite{paperinprogress}.

For the $\widehat{GL_2}$ case, we have a factorization
\[T^{-k}g_{C}^{(\alpha)}=g^{[k](\alpha)}_{-}g^{[k](\alpha)}_{+},\] where $g^{[k](\alpha)}_{-}$ is such that $\pi(g^{[k](\alpha)}_{-})=I+$ terms involving only negative powers of $z$ and $g^{[k](\alpha)}_{+}$ is such that $\pi(g^{[k](\alpha)}_{+})=A_{k}^{(\alpha)}+$ terms involving only positive powers of $z$. Here, $I$ denotes the identity matrix and $A_k^{(\alpha)}$ denotes a matrix that is independent of both $z$ and $z^{-1}$. This factorization of $T^{-k}g_{C}^{(\alpha)}$ is called the Birkhoff factorization (\cite{MR900587}). From now on, we will work in the non-centrally extended loop group, i.e. in $\widetilde{GL_2}$, but omit $\pi$ in our notation.

We define matrix Baker functions, $\Psi^{[k](\alpha)}$
by \[\Psi^{[k](\alpha)}=T^k g^{[k](\alpha)}_{-}.\]
We then have connection matrices, $U_{k}^{(\alpha)}$, between these
matrix Baker functions, given
by \[U_{k}^{(\alpha)}=(\Psi^{[k](\alpha)})^{-1}\Psi^{[k+1](\alpha)}.\]
The entries of these connection matrices have expressions in terms of
the tau-functions. (Another fact about the connections matrices, which
will be particularly important in our discussion of orthogonal
polynomials in the next section, is that they are nonnegative in $z$,
i.e. none of their entries have negative powers of
$z$. \cite{paperinprogress})

We obtain the difference relations by factoring our connection
matrices,
\[U_{k}^{(\alpha)}=V_{k}^{(\alpha)}(W_{k}^{(\alpha)})^{-1}=(W_k^{(\alpha-1)})^{-1}V_{k+1}^{(\alpha-1)},\] where \[V_{k}^{(\alpha)}=(g_{-}^{[k](\alpha)})^{-1}Q_{0}^{-1}g_{-}^{[k](\alpha+1)}\] and \[W_{k}^{(\alpha)}=(g_{-}^{[k+1](\alpha))})^{-1}Q_{1}^{-1}g_{-}^{[k](\alpha+1)}.\] The equality of these two different expressions for the connection matrices are called ``zero curvature equations" and they give us equalities satisfied by the tau-functions. These equalities imply the difference relations previously found using the Desnanot-Jacobi identity.

More explicitly, we have the following lemma, the proof of which is
included in \cite{paperinprogress}:
\begin{lemma}\cite{paperinprogress}
\[V_{k}^{(\alpha)}=\left[\begin{array}{ccc}
z-\frac{\tau_{k-1}^{(\alpha+1)}\tau_{k+1}^{(\alpha)}}{\tau_k^{(\alpha+1)}\tau_k^{(\alpha)}} & & \frac{\tau_{k-1}^{(\alpha+1)}}{\tau_{k}^{(\alpha+1)}}\\
 & & \\
-\frac{\tau_{k+1}^{(\alpha)}}{\tau_k^{(\alpha)}} & & 1
\end{array}\right], \quad W_{k}^{(\alpha)}=\left[\begin{array}{ccc}
1 & & -\frac{\tau_{k}^{(\alpha)}}{\tau_{k+1}^{(\alpha)}}\\
 & & \\
\frac{\tau_{k+1}^{(\alpha+1)}}{\tau_k^{(\alpha+1)}} & & z-\frac{\tau_{k}^{(\alpha)}\tau_{k+1}^{(\alpha+1)}}{\tau_{k+1}^{(\alpha)}\tau_k^{(\alpha+1)}}
\end{array}\right].\]
\end{lemma}
Since $V_{k}^{(\alpha)}$ and $W_{k}^{(\alpha)}$ satisfy,
\[
U_{k}^{(\alpha)}=V_{k}^{(\alpha)}(W_{k}^{(\alpha)})^{-1}=(W_{k}^{(\alpha-1)})^{-1}V_{k+1}^{(\alpha-1)},
\] 
we find that
\[\\U_{k}^{(\alpha)}=\left[\begin{array}{ccc}
z-\frac{\tau_{k}^{(\alpha)}\tau_{k+1}^{(\alpha+1)}}{\tau_{k+1}^{(\alpha)}\tau_{k}^{(\alpha+1)}}-\frac{\tau_{k-1}^{(\alpha+1)}\tau_{k+1}^{(\alpha)}}{\tau_{k}^{(\alpha+1)}\tau_{k}^{(\alpha)}} & & \frac{\tau_{k}^{(\alpha)}}{\tau_{k+1}^{(\alpha)}}\\
 & & \\
-\frac{\tau_{k+1}^{(\alpha)}}{\tau_{k}^{(\alpha)}} & & 0\\
\end{array}\right] =\left[\begin{array}{ccc}
z-\frac{\tau_{k}^{(\alpha)}\tau_{k+2}^{(\alpha-1)}}{\tau_{k+1}^{(\alpha)}\tau_{k+1}^{(\alpha-1)}}-\frac{\tau_{k}^{(\alpha-1)}\tau_{k+1}^{(\alpha)}}{\tau_{k+1}^{(\alpha-1)}\tau_{k}^{(\alpha)}} & & \frac{\tau_{k}^{(\alpha)}}{\tau_{k+1}^{(\alpha)}}\\
 & & \\
-\frac{\tau_{k+1}^{(\alpha)}}{\tau_{k}^{(\alpha)}} & 0\\
\end{array}\right].\] 
This implies
\[(\tau_{k}^{(\alpha)})^2(\tau_{k+2}^{(\alpha-1)}\tau_{k}^{(\alpha+1)}-\tau_{k+1}^{(\alpha-1)}\tau_{k+1}^{(\alpha+1)})=(\tau_{k+1}^{(\alpha)})^2(\tau_{k+1}^{(\alpha-1)}\tau_{k-1}^{(\alpha+1)}-\tau_{k}^{(\alpha-1)}\tau_{k}^{(\alpha+1)}).\]
We notice that, if
\[\tau_{k+1}^{(\alpha-1)}\tau_{k-1}^{(\alpha+1)}=\tau_{k}^{(\alpha+1)}\tau_{k}^{(\alpha-1)}-(\tau_{k}^{(\alpha)})^2\]
holds for some $k$, the above identity implies that it holds for
$k+1$. Since $\tau_{-1}^{(\alpha)}=0$ and $\tau_{0}^{\alpha-1}=1$,
$\tau_{k+1}^{(\alpha-1)}\tau_{k-1}^{(\alpha+1)}=\tau_{k}^{(\alpha+1)}\tau_{k}^{(\alpha-1)}-(\tau_{k}^{(\alpha)})^2$
holds trivially for $k=0$, and hence
\[\tau_{k+1}^{(\alpha-1)}\tau_{k-1}^{(\alpha+1)}=\tau_{k}^{(\alpha+1)}\tau_{k}^{(\alpha-1)}-(\tau_{k}^{(\alpha)})^2\]
holds for all $k\ge0$, which is exactly what we previously obtained
using the Desnanot-Jacobi identity \eqref{eq:2}.

\subsection{$\widehat{GL_3}$ Case}
The matrix Baker functions in the $\widehat{GL_3}$ case are given in a
completely analogous way to the $\widehat{GL_2}$
case. Here,
\[\Psi^{[k,\ell](\alpha,\beta)}=T_1^kT_2^{\ell}
g^{[k,\ell](\alpha,\beta)}_{-},\]
where $g^{[k,\ell](\alpha,\beta)}_{-}$ is the part of the Birkhoff
factorization of $T_{2}^{-\ell}T_{1}^{-k}g_{C,D,E}^{(\alpha,\beta)}$
that is a lift of $I$ plus a $3\times 3$ matrix whose entries have
only negative powers of $z$. (As before, $I$ denotes the identity
matrix.)

We then have two sets of connection matrices,
$U_{[k_+,\ell]}^{(\alpha,\beta)}$ and
$U_{[k,\ell_+]}^{(\alpha,\beta)},$ that we factor in two different ways
to obtain our zero curvature equations.
\[U_{[k_+,\ell]}^{(\alpha,\beta)}=(\Psi^{[k,\ell](\alpha,\beta)})^{-1}\Psi^{[k+1,\ell](\alpha,\beta)},\]
\[U_{[k,\ell_+]}^{(\alpha,\beta)}=(\Psi^{[k,\ell](\alpha,\beta)})^{-1}\Psi^{[k,\ell+1](\alpha,\beta)}.\]
Factoring these connection matrices and performing calculations
similar to those above, we find that our $\widehat{GL_3}$
tau-functions satisfy the system of four difference equations listed
in the introduction. Under suitable changes of variables, equations
$(2)$ and $(4)$ can be shown to be $T$-system relations. $T$-systems
are a sort of generalization of $Q$-systems which have
interpretations, for example in terms of perfect matchings of graphs
(\cite{MR2317336}). Like $Q$-system relations, $T$-system relations
can also be seen as cluster algebra mutations (\cite{MR2551179}). We
comment that equations $(2)$ and $(4)$ are independent of each other
($(2)$ depends only on the $k,\ell,\alpha$ parameters and equation
$(4)$ depends on the $k,\ell,\beta$ parameters,), so our tau-functions
satisfy $T$-system relations in two different ways. We believe that
the remaining two relations, $(1)$ and $(3)$, are independent of each
other and are not implied by the $T$-system relations. We are
currently working to find other situations in which equations $(1)$
and $(3)$ appear.

\section{Orthogonal Polynomials from Connection Matrices}
\subsection{$\widehat{GL_2}$ Case} In the $\widehat{GL_2}$ case, we
can use the fact that the connection matrices are nonnegative in $z$
to obtain orthogonal polynomials. In the following, we will explain
how to do this. (For more on orthogonal polynomials see, for example,
\cite{MR2542683}).

\begin{definition}\cite{paperinprogress}
  Define
  $S:\mathbb{C}[c_k]_{k\in \mathbb{Z}}\rightarrow
  \mathbb{C}[c_k]_{k\in \mathbb{Z}}$
  to be the multiplicative map such that $S(1)=0$ and $S(c_k)=c_{k+1}$
  for all $k$. The shift fields, $S^{\pm}(z)$, are then given by
  $S^{\pm}(z)=(1-\frac{S}{z})^{\pm}$, which also act multiplicatively.
\end{definition}
\begin{example}
  \[S^{+}(z)\tau_{2}^{(\alpha)}=S^{+}(z)\det\begin{bmatrix} c_{\alpha} & c_{\alpha+1}\\
    c_{\alpha+1} &
    c_{\alpha+2}\end{bmatrix}=(c_{\alpha}-c_{\alpha+1}/z)(c_{\alpha+2}-c_{\alpha+3}/z)-(c_{\alpha+1}-c_{\alpha+2}/z)^2=\]
\[=\det\begin{bmatrix} c_{\alpha} & c_{\alpha+1}\\
c_{\alpha+1} & c_{\alpha+2}\end{bmatrix}-\det\begin{bmatrix} c_{\alpha} & c_{\alpha+1}\\
c_{\alpha+2} & c_{\alpha+3}\end{bmatrix}/z+\det\begin{bmatrix} c_{\alpha+1} & c_{\alpha+2}\\
c_{\alpha+2} & c_{\alpha+3}\end{bmatrix}/z^{2}=\]
\[=\tau_{2}^{(\alpha)}-\det\begin{bmatrix} c_{\alpha} & c_{\alpha+1}\\
c_{\alpha+2} & c_{\alpha+3}\end{bmatrix}/z+\tau_{2}^{(\alpha+1)}/z^{2}\]
\end{example}
We note that the positive shift field, $S^{+}(z)$ sends elements of
$\mathbb{C}[c_{k}]_{k\in\mathbb{Z}}$ to polynomials in $z^{-1}$ with
coeffients in $\mathbb{C}[c_k]_{k\in\mathbb{Z}}$. In particular, as is
illustrated in the above example, $S^{+}(z)$ sends
$\tau_{k}^{(\alpha)}$ to a Laurent polynomial in $z$, with smallest
degree equal to $-k$ and largest degree equal to $0$. We also have the
following useful formula:

\begin{equation}
\label{eq:7}
z^{k}S^{+}(z)\tau_{k}^{(\alpha)}=\det \begin{bmatrix}
c_{\alpha} & c_{\alpha+1}&\cdots &c_{\alpha+k-1} & 1\\
c_{\alpha+1} & c_{\alpha+2}&\cdots &c_{\alpha+k} & z\\
\vdots & \vdots & \cdots & \vdots & \vdots\\
c_{\alpha+k} & c_{\alpha+k+1} & \cdots & c_{\alpha+2k-1}& z^{k}
\end{bmatrix}\end{equation} for all $\alpha\in \mathbb{Z}$ and for all
$k\ge 0$. (This comes from the fact
that
\[z^{k}S^{+}(z)\tau_{k}^{(\alpha)}=\frac{1}{k!}\Res_{\mathbf{w}}(\displaystyle
\prod_{1\le i<j\le k}(w_i-w_j)^2\displaystyle
\prod_{i=1}^{k}(z-w_i)\prod_{i=1}^{k}C^{(\alpha)}(w_i)).)\]
Observe that formulas like \eqref{eq:7} appear in the theory of
orthogonal polynomials. (See for example, equation $(2.1.6)$ in
\cite{MR2542683}). Below, we will use connection matrices to derive
the orthogonality of the polynomials given by \eqref{eq:7}.

The negative shift field, $S^{-}(z)$, sends elements of  $\mathbb{C}[c_{k}]_{k\in\mathbb{Z}}$ to series in $z^{-1}$ with coeffients in  $\mathbb{C}[c_{k}]_{k\in\mathbb{Z}}$. For example, $S^{-}(z)c_{\alpha}=\displaystyle\sum_{i=0}^{\infty}\frac{c_{\alpha+i}}{z^{i}}$.

\begin{theorem}\cite{paperinprogress}
\label{thm:4}
$\pi(g_{-}^{[k],(\alpha)})=\frac{1}{\tau_{k}^{(\alpha)}}\begin{bmatrix}
S^{+}(z)\tau_{k}^{(\alpha)} & S^{+}(z)\tau_{k-1}^{(\alpha)}/z\\
S^{-}(z)\tau_{k+1}^{(\alpha)}/z & S^{-}(z)\tau_{k}^{(\alpha)}
\end{bmatrix}$
\end{theorem}

Recall the definition of the Baker functions, $\Psi^{[k](\alpha)}=T^{k}g_{-}^{[k](\alpha)}$. Any two Baker functions are related by elementary connection matrices. In particular, \[(\Psi^{[0](\alpha)})^{-1}\Psi^{[k](\alpha)}=U_{0}^{(\alpha)}U_{1}^{(\alpha)}\cdots U_{k-1}^{(\alpha)}\] is a product of connection matrices and so is nonnegative in $z$.

Using the above theorem, we see that 
\begin{equation*}
(\Psi^{[0](\alpha)})^{-1}\Psi^{[k](\alpha)}=\begin{bmatrix}1 & 0\\
-\displaystyle \sum_{i=0}^{\infty}\frac{c_{\alpha+i}}{z^{i+1}} & 1
\end{bmatrix} \frac{1}{\tau_{k}^{(\alpha)}}
\begin{bmatrix}
z^{k}S^{+}(z)\tau_{k}^{(\alpha)} & z^{k-1}S^{+}(z)\tau_{k-1}^{(\alpha)}\\
z^{-k-1}S^{-}(z)\tau_{k+1}^{(\alpha)} & z^{-k}S^{-}(z)\tau_{k}^{(\alpha)}
\end{bmatrix}.
\end{equation*}
Given $\alpha$, denote by $\langle\,,\,\rangle$ the bilinear
product given by
\begin{equation}
 \label{eq:11}\langle
  f(z),g(z)\rangle=\Res_{z}(\displaystyle
  \sum_{i=0}^{\infty}\frac{c_{\alpha+i}}{z^{i+1}}f(z)g(z)),\end{equation}
for all polynomials $f(z)$ and $g(z)$, and
denote 
\begin{equation}
\label{eq:4}
p_{k}^{(\alpha)}(z)=\frac{1}{\tau_{k}^{(\alpha)}}z^{k}S^{+}(z)\tau_{k}^{(\alpha)},
\end{equation}
which we note is a monic polynomial of degree $k.$ Consider the entry
in the first column and second row of
$(\Psi^{[0](\alpha)})^{-1}\Psi^{[k](\alpha)}$. Since
$(\Psi^{[0](\alpha)})^{-1}\Psi^{[k](\alpha)}$ is nonnegative in $z$
and $z^{-k-1}S^{-}(z)\tau_{k+1}^{(\alpha)}$ has highest degree equal
to $-k-1,$ we see that 
\[\langle p_{k}^{(\alpha)}(z),z^{n}\rangle=0
\]
for all $0\le n<k$. So
\[\langle p_{k}^{(\alpha)}(z),p_{\ell}^{(\alpha)}(z)\rangle=0
\]
for $k \ne \ell$ and the $\{p_{k}^{(\alpha)}\}$ are a collection of
orthogonal polynomials. Since the series $C^{(\alpha)}(z)$ can be
defined arbitrarily, we in fact obtain all orthogonal polynomials in
this way. So the theory of orthogonal polynomials appears as a subset
in the study of the representation theory of $\widehat{GL_2}$.
\begin{example}
  Hermite polynomials are the monic polynomials orthogonal for the
  following bilinear form:
\begin{equation}\label{eq:8}
\langle f(z),g(z)\rangle=\int_{-\infty}^{\infty}f(z)g(z)e^{z^{-2}}dz,
\end{equation} where the integral is along the real axis. The moments
for this bilinear form are given by 
\begin{equation}
c_i=\int_{-\infty}^{\infty}z^ie^{-z^2}dz=\begin{cases}
0 &i\text{ is odd}\\
\frac{(2m)!\sqrt{\pi}}{m!2^m} & i=2m.
\end{cases} 
\end{equation} We can then rewrite \eqref{eq:8} as
\begin{equation}
\langle f(z),g(z)\rangle=\Res_{z}(\sum_{i=0}^{\infty}\frac{c_i}{z^{i+1}}f(z)g(z)),
\end{equation} 
which is precisely \eqref{eq:11} for our specified moments, $c_i$,
$i\ge 0$, for $\alpha=0$. So, in this case, the orthogonal polynomials given by (\ref{eq:4}) are
precisely the Hermite polynomials (see, \cite{MR2542683}).
\end{example}
\subsection{$\widehat{GL_3}$ Case}
We obtain orthogonal polynomials from the $\widehat{GL_3}$ case
exactly as we did in the $\widehat{GL_2}$ case, by using the
connection matrices. For this, we need the analogue of Theorem
\ref{thm:4} which is as follows. (This theorem is stated and
proven in \cite{paperinprogress}.):
\begin{theorem}
$\pi(g_{-}^{[k,\ell](\alpha,\beta)})=
\\
\frac{1}{\tau_{k,\ell}^{(\alpha,\beta)}}\begin{bmatrix}
S_{c}^{+}(z)S_d^{+}(z)\tau_{k,\ell}^{(\alpha,\beta)} & \frac{S_{c}^{+}(z)S_d^{+}(z)\tau_{k-1,\ell}^{(\alpha,\beta)}}{z} & (-1)^k \frac{S_{c}^{+}(z)S_d^{+}(z)\tau_{k-1,\ell-1}^{(\alpha,\beta)}}{z}\\
\frac{S_{c}^{-}(z)S_e^{+}(z)\tau_{k+1,\ell}^{(\alpha,\beta)}}{z} & S_{c}^{-}(z)S_e^{+}(z)\tau_{k,\ell}^{(\alpha,\beta)} & \frac{S_{c}^{-}(z)S_e^{+}(z)\tau_{k,\ell-1}^{(\alpha,\beta)}}{z}\\
(-1)^{k+1}\frac{S_{d}^{-}(z)S_e^{-}(z)\tau_{k+1,\ell+1}^{(\alpha,\beta)}}{z} & (-1)^k \frac{S_{d}^{-}(z)S_e^{-}(z)\tau_{k,\ell+1}^{(\alpha,\beta)}}{z} & S_{d}^{-}(z)S_e^{-}(z)\tau_{k,\ell}^{(\alpha,\beta)} \\
\end{bmatrix}$
\end{theorem} Here, the shift operators, $S_{c}^{\pm}(z)$, $S_d^{\pm}(z)$, and $S_e^{\pm}(z)$ are analagous to the shift operators defined in the $\widehat{GL_2}$ case. $S_{c}^{\pm}(z)$ acts on the $c_i$s exactly as $S^{\pm}(z)$ does in the $\widehat{GL_2}$ case and acts trivially on the $d_i$s and $e_i$s. $S_d^{\pm}(z)$ and $S_e^{\pm}(z)$ are similarly defined.

The matrix Baker functions for the $\widehat{GL_3}$ case are given above, by $\Psi^{[k,\ell](\alpha,\beta)}=T_1^kT_2^{\ell} g^{[k,\ell](\alpha,\beta)}_{-}$. We then see that
\[(\Psi^{[0,0](\alpha,\beta)})^{-1}\Psi^{[k,\ell](\alpha,\beta)}\] is nonnegative in $z$, since it is equal to a product of connection matrices: \[(\Psi^{[0,0](\alpha,\beta)})^{-1}\Psi^{[k,\ell](\alpha,\beta)}=U_{[0_+,0]}^{(\alpha,\beta)}U_{[1_+,0]}^{(\alpha,\beta)}\cdots U_{[{k-1}_+,0]}^{(\alpha,\beta)}U_{[k,0_+]}^{(\alpha,\beta)}U_{[k,1_+]}^{(\alpha,\beta)}\cdots U_{[k,{\ell-1}_+]}^{(\alpha,\beta)}.\]

In the $\widehat{GL_2}$ case, our orthogonal polynomials were obtained
by acting on tau-functions by the shift operators. Since our
$\widehat{GL_3}$ tau-functions are, in general, infinite sums, to
obtain polynomials we restrict to the case that either the series
$C^{(\alpha-\beta)}(z)$ or the series $E^{(\beta)}(z)$ is $0$. Here, we discuss the case that
$E^{(\beta)}(z)=0$.

Using the formula (in terms of residues) for the $\widehat{GL_3}$
tau-functions given earlier, we can show that when
$E^{(\beta)}(z)=0$, \[\tau_{k,\ell}^{(\alpha,\beta)}=0 \text{ when }k<\ell \]
and when $k\ge
\ell$,
\[\tau_{k,\ell}^{(\alpha,\beta)}=(-1)^{\frac{\ell(\ell+1)}{2}}\det\begin{bmatrix}
  d_{\alpha} & \cdots &d_{\alpha+\ell-1}&c_{\alpha-\beta} &\cdots &c_{\alpha-\beta+k-\ell-1}\\
  d_{\alpha+1} &\cdots &d_{\alpha+\ell}&c_{\alpha-\beta+1} &\cdots &c_{\alpha-\beta+k-\ell}\\
  \vdots & \cdots &\vdots & \vdots & \cdots & \vdots\\
  d_{\alpha+k-1} & \cdots &d_{\alpha+k-\ell-2}&c_{\alpha-\beta+k-1}
  &\cdots &c_{\alpha-\beta+2k-\ell-2}\end{bmatrix}.\]

To see this, note that when $E^{(\beta)}(z)=0$, $p_{n_c,n_d,n_e}$ in the formula
for our $\tau_{k,\ell}^{(\alpha,\beta)}$ reduces to
\[p_{n_c,n_d}=\displaystyle (-1)^{\displaystyle\frac{n_d(n_d+1)}{2}}\times\]
\[\displaystyle \prod_{1\le i<j\le n_c}(x_i-x_j)^2 \displaystyle
\prod_{1\le i<j\le n_d}(y_i-y_j)^2
\prod_{i=1}^{n_c}\prod_{j=1}^{n_d}(x_i-y_j),\]
which is a Vandermonde determinant in the $y_i$s and $x_i$s
times
\[\prod_{1\le i<j\le n_c}(x_i-x_j) \prod_{1\le i<j\le n_d}(y_i-y_j).\]

The formula for these tau functions is no longer a sum of
$c_{n_c,n_d,n_e}$s, since $n_e$ is zero and, when $k\ge\ell$, there is
only one choice for $n_c$ and $n_d$ such that $n_c+n_d=k$ and
$n_d=\ell$. When $k<\ell$, no choice of $n_c$ and $n_d$ exists. We
then have, when $k\ge\ell$,
\[\tau_{k,\ell}^{(\alpha,\beta)}=\frac{1}{(k-\ell)!\ell!}\Res_{\mathbf{x}}\Res_{\mathbf{y}}\left(\displaystyle\prod_{i=1}^{k-\ell}C^{(\alpha-\beta)}(x_{i})\displaystyle\prod_{i=1}^{\ell}D^{(\alpha)}(y_{i})p_{k-\ell,\ell}\right).\]
Since we have no $e_i$s in our formulas, the $S^{\pm}_e(z)$ act
trivially on our tau functions, so we only need $S_{c}^{\pm}(z)$ and
$S_{d}^{\pm}(z)$. So we have
$z^{k}S_{c}^{+}(z)S_{d}^{+}(z)\tau_{k,\ell}^{(\alpha,\beta)}=$\[(-1)^{\frac{\ell(\ell+1)}{2}}\det
\begin{bmatrix}
  d_{\alpha} & \cdots &d_{\alpha+\ell-1}&c_{\alpha-\beta} &\cdots &c_{\alpha-\beta+k-\ell-1} & 1\\
  d_{\alpha+1} &\cdots &d_{\alpha+\ell}&c_{\alpha-\beta+1} &\cdots &c_{\alpha-\beta+k-\ell} & z \\
  \vdots & \cdots &\vdots & \vdots & \cdots & \vdots & \vdots\\
  d_{\alpha+k} & \cdots &d_{\alpha+k-\ell-1}&c_{\alpha-\beta+k}
  &\cdots &c_{\alpha-\beta+2k-\ell-1}& z^{k}\end{bmatrix},\]
which comes from the fact that
\[S^{+}_c(z)
S^{+}_d(z)\tau_{k,\ell}^{(\alpha,\beta)}=\frac{1}{(k-\ell)!\ell!}\Res_{\mathbf{x}}\Res_{\mathbf{y}}\left(\displaystyle\prod_{i=1}^{k-\ell}C^{(\alpha-\beta)}(x_{i})\displaystyle\prod_{i=1}^{\ell}D^{(\alpha)}(y_{i})\prod_{i=1}^{k-\ell}(z-x_i)\prod_{i=1}^{\ell}(z-y_i)p_{k-\ell,\ell}\right)\]
and
\[\prod_{i=1}^{k-\ell}(z-x_i)\prod_{i=1}^{\ell}(z-y_i)p_{k-\ell,\ell}\]
is a Vandermonde determinant in the $y_i$s, $x_i$s, and $z$ times
\[\prod_{1\le i<j\le k-\ell}(x_i-x_j)\prod_{1\le i<j\le
  \ell}(y_i-y_j).\]

Denote
$p_{k,\ell}^{(\alpha,\beta)}(z)=\frac{1}{\tau_{k,\ell}^{(\alpha,\beta)}}z^{k}S_{c}^{+}(z)S_{d}^{+}(z)\tau_{k,\ell}^{(\alpha,\beta)}$,
which is a monic polynomial of degree $k$.
Given $\alpha,\beta\in \mathbb{Z}$,
 define $\langle\,, \,\rangle_{C}$ to be the bilinear product given
 by \[\langle f(z), g(z)\rangle_{C}=\Res_{z}(\displaystyle
 \sum_{i=0}^{\infty}\frac{c_{\alpha-\beta+i}}{z^{i+1}}f(z)g(z))=0,\]
 for all polynomials $f(z)$ and $g(z)$. Similarly, define
 $\langle\,,\,\rangle_{D}$ to be the bilinear product given by 
\[
\langle f(z), g(z)\rangle_{D}=\Res_{z}(\displaystyle
\sum_{i=0}^{\infty}\frac{d_{\alpha+i}}{z^{i+1}}f(z)g(z))=0.
\] 
Consider the second row and first column of 
\begin{multline*}
  (\Psi^{[0,0](\alpha,\beta)})^{-1}\Psi^{[k,\ell](\alpha,\beta)}=
  \begin{bmatrix}1 & 0 & 0\\
    -\displaystyle\sum_{i=0}^{\infty}\frac{c_{\alpha-\beta+i}}{z^{i+1}}& 1 & 0\\
    -\displaystyle \sum_{i=0}^{\infty}\frac{d_{\alpha+i}}{z^{i+1}} & 0 & 1\\
  \end{bmatrix}\times
  \\
  \frac{1}{\tau_{k,\ell}^{(\alpha,\beta)}}
  \begin{bmatrix}
    z^kS_{c}^{+}(z)S_d^{+}(z)\tau_{k,\ell}^{(\alpha,\beta)} &
    z^k\frac{S_{c}^{+}(z)S_d^{+}(z)\tau_{k-1,\ell}^{(\alpha,\beta)}}{z}
    &
    (-1)^k z^k\frac{S_{c}^{+}(z)S_d^{+}(z)\tau_{k-1,\ell-1}^{(\alpha,\beta)}}{z}\\
    z^{\ell-k}\frac{S_{c}^{-}(z)\tau_{k+1,\ell}^{(\alpha,\beta)}}{z} &
    z^{\ell-k}S_{c}^{-}(z)\tau_{k,\ell}^{(\alpha,\beta)} &
    z^{\ell-k}\frac{S_{c}^{-}(z)\tau_{k,\ell-1}^{(\alpha,\beta)}}{z}\\
    (-1)^{k+1}z^{-\ell}\frac{S_{d}^{-}(z)\tau_{k+1,\ell+1}^{(\alpha,\beta)}}{z}
    & (-1)^k
    z^{-\ell}\frac{S_{d}^{-}(z)\tau_{k,\ell+1}^{(\alpha,\beta)}}{z}
    & z^{-\ell}S_{d}^{-}(z)\tau_{k,\ell}^{(\alpha,\beta)} \\
  \end{bmatrix}.
\end{multline*}
Since
$z^{\ell-k}\frac{S_{c}^{-}(z)\tau_{k+1,\ell}^{(\alpha,\beta)}}{z}$ has
highest degree $\ell-k-1$, we see that 
\[\langle p_{k,\ell}^{(\alpha,\beta)},z^{n}\rangle_{C}=0\]
for $0\le n < k-\ell$. Similarly, the entry in the third row and first
column gives
 \[\langle p_{k,\ell}^{(\alpha,\beta)},z^{n}\rangle_D=0
\]
for $0\le n < \ell$, since
$(-1)^{k+1}z^{-\ell}\frac{S_{d}^{-}(z)\tau_{k+1,\ell+1}^{(\alpha,\beta)}}{z}$
has degree at most $-\ell-1$. Such polynomials,
$p^{(\alpha,\beta)}_{k,\ell}(z),$ are known as ``type II multiple
orthogonal polynomials'' (see \cite{MR2542683}, \cite{MR2006283}).

We plan to continue studying orthogonal polynomials coming from the $\widehat{GL_3}$ case. In particular, we would like to better understand more general cases, in which the series $E^{(\beta)}(z)$ is not required to be $0$.
\def\cprime{$'$}
\providecommand{\bysame}{\leavevmode\hbox to3em{\hrulefill}\thinspace}
\providecommand{\MR}{\relax\ifhmode\unskip\space\fi MR }

\providecommand{\MRhref}[2]{
  \href{http://www.ams.org/mathscinet-getitem?mr=#1}{#2}
}
\providecommand{\href}[2]{#2}

\end{document}